\documentclass[12pt,a4paper,leqno]{article}
\usepackage {amsmath, latexsym, amscd, amssymb}
\usepackage {graphicx}
\usepackage[latin1]{inputenc}

\def\tvi #1#2{\vrule height #1pt depth #2pt width 0mm}
\def\build#1_#2^#3{\mathrel{\mathop{\kern 0pt#1}\limits_{#2}^{#3}}}

\newcommand{\C}{{\mathbb{C}}}

\newcommand{\R}{{\mathbb{R}}}

\newcommand{\sB}{\mathcal{B}}

\newcommand{\sS}{\mathcal{S}}
\newcommand{\sL}{\mathcal{L}}

\newcommand{\sO}{\mathcal{O}}

\newcommand{\hf}{\widehat f}
\newcommand{\hF}{\widehat F}
\newcommand{\hg}{\widehat g}

\newcommand{\hu}{\widehat u}
\newcommand{\hv}{\widehat v}
\newcommand{\hw}{\widehat w}

\newcommand{\Ot}{\sO(D_\rho)\lbrack\lbrack t\rbrack\rbrack}
\newcommand{\Ots}{\sO(D_\rho)\lbrace\lbrace t\rbrace\rbrace}

\newcommand{\cf}{{\it cf.} }
\newcommand{\ie}{{\it i.e.}, }

\newcommand{\proof}{ \noindent { \em {\sc Proof.} }}
\newcommand{\qed}{ $\hfill {\Box}$\break}

\newcommand{\gvf }{\varphi }

\newcommand{\gz}{\zeta}

\newtheorem{lemma}{Lemma}[section]
\newtheorem{thm}[lemma]{Theorem}
\newtheorem{cor}[lemma]{Corollary}
\newtheorem{pro}[lemma]{Proposition}
\newtheorem{defi}[lemma]{Definition}

\newtheorem{cexa}[lemma]{Counter-example}

\makeindex
\title{Summability of solutions of the heat equation with
inhomogeneous thermal conductivity in two variables}
\author{Werner BALSER\\ Abteilung Angewandte Analysis\\
Universit\"at Ulm, D-89069 ULM, Germany\\
{\it Email}: balser@mathematik.uni-ulm.de\\\\
Mich\`ele LODAY-RICHAUD\\ LAREMA, Universit\'e d'Angers,
2 boulevard Lavoisier\\ 49 045 ANGERS cedex 01, France
\\{\it Email}: michele.loday@univ-angers.fr}

\begin{document}

\baselineskip=1.6em
\maketitle
%\null\vspace{3cm}

\abstract{We investigate Gevrey order and 1-summability properties of the formal solution of a general heat equation in two variables. In particular, we give necessary and sufficient conditions for the 1-summability of the solution in a given direction. When restricted to the case of constants coefficients, these conditions coincide with those given by D.A. Lutz, M. Miyake, R. Sch\"afke in a 1999 article (\cite{LMS99}), and we thus provide a new proof of their result.
}
\tableofcontents
%\include{bibliographie}
%\printindex

\baselineskip=1.6em

\section{The problem}
A formal solution of the classical heat initial conditions
problem
\begin{equation}\label{classical}
 \left\lbrace
 \begin{array}{l}
\partial_t\, u-\partial^2_z\, u=0\\
\noalign{\medskip}
u(0,z)=\gvf(z)\\
\end{array}
\right.
\end{equation}
in one dimensional spatial variable $z$ reads in the form
\begin{equation*}
\begin{array}{lcl}
\hu(t,z)&=&\displaystyle
\exp\big(t\,\partial_z^2t\big)\gvf(z)\\
\noalign{\medskip}
&=&\displaystyle \sum_{j\geq 0}\frac{t^j}{j!}\gvf^{(2j)}(z)
\\
\end{array}
\end{equation*}
provided that all derivatives $\gvf^{(2j)}$
exist\footnote{We denote $\hu$, with a hat, to emphasize the
possible divergence of the series
$\hu$.}. When $\gvf\in\sO(D_p)$ is holomorphic
in a disc $D_\rho$ with center 0 and radius $\rho$ and hence
satisfies, for any $r<\rho$,  estimates of the form
\begin{equation*}
\big\vert \gvf^{(2j)}(z) \big\vert \leq C\, K^{2j}\, \Gamma(1+2j)! ,
\end{equation*} for all $j\geq 0$ and positive constants
$C$ and
$K$, on $D_r$ then, $\hu(t,z)\in\sO(D_\rho)\lbrack\lbrack t\rbrack\rbrack$
is a series of Gevrey type of order 1 in $t$ for all $z\in D_\rho$ (in
short, a 1-Gevrey series). The
Gevrey estimates are locally uniform with respect to $z$ in
$D_\rho$. These conditions are optimal as shown by the following
example: Let consider  $\displaystyle
 \tvi{18}{15}\gvf(z)=\frac{1}{1-z}=\sum_{n\geq 0}z^n$ so
that
 $\gvf^{(2j)}(0)=(2j)!$. The corresponding solution
$\hu(t,z)$ is  of exact Gevrey order 1 and,
in particular, is divergent.  It turns out
that it is actually 1-summable in
all  direction but $\R^+$ in the sense of Definition
\ref{1sum} below, that is, 1-summable in $t$  uniformally
with respect to $z$ near 0.

In 1999, D. Lutz, M. Miyake and R. Sch\"afke (\cite{LMS99}) gave
necessary and sufficient conditions on $\gvf$ for $\hu$ to
be 1-summable in a given direction $\arg t=\theta$. Various works have been done towards the summability of divergent solutions of  partial differential equations with constant coefficients ( \cite{Bal99}, \cite{Miy99}, \cite{BM99}, \cite{Bal04},\dots) or variable coefficients (\cite{H99}, \cite{Ou02}, \cite{PZ97}, \cite{Mk08}, \cite{Mk09},\dots) in two variables. In \cite{Mk05}, S. Malek has investigated the case of linear partial differential equations with constant coefficients in more variables.
\bigskip

In this article we are interested in the very general heat initial
conditions problem with inhomogeneous thermal conductivity
 and internal heat generation
\begin{equation}\label{inhomheatq}
\left\lbrace\begin{array}{l}
\partial_t\, u-a(z)\,\partial_z^2\, u=q(t,z)\qquad
a(z)\in\sO(D_\rho)\\
\noalign{\medskip}
u(0,z)=\gvf(z)\in\sO(D_\rho).\\
\end{array}\right.
\end{equation}
The heat equation describes heat propagation under
thermodynamics and Fourier laws. The coefficient
$a(z)$, named thermal diffusivity, is related to the
thermal conductivity
$\kappa$ by the formula
$\displaystyle a=\frac{\kappa}{c\rho}$ where $c$ is the capacity and
$\rho$ the density of the medium. We assume that $a(z)$ and $\gvf(z)$ are analytic on
a neighborhood of $z=0$.
The internal
heat input $q$ may be smooth or not. An important case is the case with no internal heat
generation corresponding to a homogeneous heat equation:
\begin{equation}\label{homheat}
\left\lbrace\begin{array}{l}
\partial_t\, u-a(z)\,\partial_z^2\, u=0\qquad
a(z)\in\sO(D_\rho)\\
\noalign{\medskip}
u(0,z)=\gvf(z)\in\sO(D_\rho).\\
\end{array}\right.
\end{equation}
 In case of an isotropic and
homogeneous medium, $\kappa, c, \rho$ and hence $a$ are constants. An
adequate choice of units allows then to assume $a=1$ and the equation
reduces to the reference heat equation $\partial_t u -\partial_z^2 u=0$.
\bigskip

Actually,   for
notational convenience, we consider the
 problem in the form
\begin{equation}\label{Dheat}
\begin{tabular}
{| p{5.8cm}  |}
\hline  $\displaystyle
\quad
\big(1-a(z)\,\partial_t^{-1}\partial_z^2\big)\,\hu=\hf(t,z)$
\tvi{20}{15}\\ \hline
\end{tabular} \ ,\
a(z)\in\sO(D_\rho)\text{ and }\hf(t,z)\in
\sO(D_\rho)\lbrack\lbrack t\rbrack\rbrack
\end{equation}
where $\partial_t^{-1}\,\hu$ stands for the anti-derivative
$\int_0^t \hu(s,z)ds$ of $\hu$ with respect to $t$ which
vanishes at $t=0$.

 Problem (\ref{Dheat}) is
equivalent to
\begin{equation*}
\left\lbrace\begin{array}{l}
\partial_t\, \hu-a(z)\,\partial_z^2\hu=\partial_t\,\hf(t,z)\\
\noalign{\medskip}
\hu(0,z)=\hf(0,z).\\
\end{array}\right.
\end{equation*}
and hence to Problem (\ref{inhomheatq}) by choosing
$q(t,z)=\partial_t\hf(t,z)$ and $\gvf(z)=\hf(0,z)$.\\
Moreover,  Problem (\ref{Dheat}) reduces to the homogeneous
case (\ref{homheat}) if and only if the inhomogenuity
$\hf$ does not depend on
$t$.
\bigskip

From now, we denote $D=1-a(z)\,\partial_t^{-1}\partial_z^2$
and, given a series $\hu\in\sO(D_\rho)\lbrack\lbrack
t\rbrack\rbrack$, we denote
$$
\hu(t,z)=\sum_{j\geq 0}\frac{t^j}{j!}u_{j,*}(z)=
\sum_{n\geq 0}\hu_{*,n}(t)\frac{z^n}{n!}=\sum_{j,n\geq 0}
\hu_{j,n}\frac{t^j}{j!}\frac{z^n}{n!}\cdot$$

Since $\big(\sO(D_\rho)\lbrack\lbrack
t\rbrack\rbrack,\partial_t,\partial_z\big)$  is a
differential algebra and $a(z)\in\sO(D_\rho)$ the operator
$D$ acts inside
$\sO(D_\rho)\lbrack\lbrack t\rbrack\rbrack$. More precisely, we can state:

\begin{pro}\label{linisof} The map
$$
D :\ \sO(D_\rho)\lbrack\lbrack t\rbrack\rbrack
\longrightarrow \sO(D_\rho)\lbrack\lbrack t\rbrack\rbrack
$$
is a linear isomorphism.
\end{pro}
\proof
The operator $D$ is linear. A series $\displaystyle
\hu(t,z)=\sum_{j\geq 0}\frac{t^j}{j!}\hu_{j,*}(z)$ is a
solution of Problem (\ref{Dheat}) if and only if
\begin{equation}\label{recurrence}
%\left\lbrace\begin{array}{l}
\hu_{j,*}(z)=\hf_{j,*}(z)+a(z)\, \hu_{j-1,*}''(z) \quad
\text{for all } j\geq 0 \text{ starting from } \hu_{-1,*}(z)\equiv 0.
\end{equation}
Consequently, to any $\hf(t,z)\in\sO(D_\rho\lbrack\lbrack
t\rbrack\rbrack$ there is a unique solution
$\hu(t,z)\in\sO(D_\rho\lbrack\lbrack t\rbrack\rbrack$, which
proves that $D$ is bijective.\qed

In Section 2 we show that the inhomogenuity $\hf(t,z)$ and the unique
solution $\hu(t,z)$ are together 1-Gevrey.

In Section 3 we prove necessary and sufficient conditions
for $\hu$ to be 1-summable in a given direction $\arg
t=\theta$. The conditions are valid in the case when
either $a(0)\neq 0$ or $a'(0)\neq 0$.
When $a(z)=O(z^2)$ an easy counter-example shows that even
the rationality of $\hf(t,z)$ is insufficient.

In Section 4 we discuss the accessibility of our necessary
and sufficient conditions. Indeed, the conditions are
given not only in terms of the data $\hf$ but also
in terms of the first two terms $\hu_{*,0}$ and $\hu_{*,1}$ of
the solution $\hu$ itself.\\
In the particular case $a=1$ our conditions  coincide with those of
\cite{LMS99}. We thus provide a new proof of the result of
\cite{LMS99}.
\section{Gevrey properties}

In this article, we consider  $t$ as the variable
and $z$ as a parameter. The classical notion of a series of
Gevrey type of order 1 is extended to $z$-families  as
follows.
\begin{defi}[1-Gevrey series]\label{Gevrey1}
A series $\displaystyle \hu(t,z)=\sum_{j\geq 0}\frac{t^j}{j!}\hu_{j,*}(z)\in
\Ot$ is of {\rm Gevrey type of order 1} if there exist $0<r\leq
\rho,\, C>0,\, K>0 $ such that for all $j\geq 0$ and $\vert
z\vert\leq r$ we have
$$
\vert \hu_{j,*}(z)\vert \leq C\,K^j\,\Gamma(1+2j).
$$
\end{defi}
In other words, $\hu(t,z)$ is 1-Gevrey in $t$, uniformally in
$z$ on a neighbourhood of $z=0$.

We denote $\Ot_1$ the subset of $\Ot$ made of the series
which are of Gevrey type of order 1.

\begin{pro}\label{diffal}
$\big( \Ot_1, \partial_t,\partial_z\big)$ is a differential
algebra stable under $\partial_t^{-1}$ and $\partial_z^{-1}$.
\end{pro}

\proof The proof is similar to the one without parameter.
Stability under $\partial_z$ is proved using the Cauchy
Integral Formula and is guaranted by the condition ``there
exist
$r\leq \rho\ \dots$'' in Definition \ref{Gevrey1}.\qed

It results from this Proposition that the operator
$D=1-a(z)\partial_t^{-1}\partial_z^{2}$ acts inside the
space $\Ot_1$.\\

Because  the main result of this section (Theorem \ref{isomGevrey})
 is set up using Nagumo norms
on $\sO(D_\rho)$ we begin with a recall of their definition
and main properties and we refer to \cite{N42} or to \cite{CRSS00}
for more details.

\begin{defi}[Nagumo norms]\quad
\\
Let $f\in \sO(D_\rho)$, \ $p\geq 0, \ 0<r\leq \rho$ and let
$d_r(z)=\vert z\vert- r$ denote the euclidian distance of
$z$ to the boundary of the disc $D_r$.

The Nagumo norm $\Vert f\Vert_{p,r}$ of $f$ is defined by
$$
\begin{tabular}
{| p{5cm}  |}
\hline  $\displaystyle
\quad \Vert f\Vert_{p,r}=\sup_{\vert z\vert <r}\big\vert
f(z)d_r(z)^p\big\vert$
\tvi{20}{20}\\ \hline
\end{tabular} \ .
$$
\end{defi}

\begin{pro}[Properties of Nagumo norms]\label{propNagumo}\quad

\begin{description}
\item{(i)} $\Vert  .\Vert_{p,r}$ is a norm on $\sO(D_\rho)$;
\item{(ii)} For all $z\in D_r, \ \ \vert f(z)\vert\leq \Vert
f\Vert_{p,r} d(z)^{-p}$;
\item{(iii)} $\displaystyle \Vert f\Vert_{0,r}=sup_{z\in D_r}\vert
f(z)\vert$ is the usual sup-norm on $D_r$;
\item{(iv)} $\Vert fg\Vert_{p+q,r}\leq \Vert f\Vert_{p,r} \Vert
g\Vert_{q,r}$;
\item{(v)} (most important) $\Vert f'\Vert_{p+1,r}\leq e(p+1)\Vert
f\Vert_{p,r}$.
\end{description}
\end{pro}
Note that the same index $r$ occurs on both sides of the inequality
{\it (v)}. One
gets thus an estimate of the derivative $f'$ in terms of $f$
without having to shrink the domain $D_r$.

\begin{thm}\label{isomGevrey}
The map
\begin{equation*}
D:\left\lbrace\begin{array}{ccl}
 \Ot_1 &\longrightarrow& \Ot_1\\
 \noalign{\smallskip}
\hu(t,z)&\mapsto&\hf(t,z)=D\hu(t,z)\\
\end{array}\right.
\end{equation*}
is a linear isomorphism.
\end{thm}
\proof
It results from Proposition \ref{diffal} that
$D\big(\Ot_1\big)\subset \Ot_1$ and from Proposition
\ref{linisof} that $D$ is linear and injective. We are left to prove that $D$
is also surjective.

Let $\displaystyle \hf(t,z)=\sum_{j\geq 0}
\frac{t^j}{j!}\hf_{j,*}(z)\in
\Ot_1$. The coefficients $\hf_{j,*}(z)$ satisfy
\begin{equation*}
\left\lbrace\begin{array}{l}
\bullet \hf_{j,*}(z)\in\sO(D_\rho) \text{ for all } j\geq 0.\\
\noalign{\medskip}
\bullet \text{ There exist } 0<r\leq \rho,\, C>0,\, K>0 \text{
such that for all } j\geq 0 \text{ and } \vert z\vert\leq r\\
\noalign{\medskip}
\quad \vert \hf_{j,*}(z)\vert\leq CK^j\Gamma(1+2j)!
\\
\end{array}\right.
\end{equation*}
and we look forward to similar conditions on the
coefficients $\displaystyle \hu_{j,*}(z)$ of
$\displaystyle \hu(t,z)=\sum_{j\geq 0}\frac{t^j}{j!}\, \hu_{j,*}(z)$.

From the recurrence relation (\ref{recurrence}) the relation
$$
\frac{\hu_{j,*}(z)}{\Gamma(1+2j)}=\frac{\hf_{j,*}(z)}{\Gamma(1+2j)} +
a(z)\frac{\hu''_{j-1,*}(z)}{\Gamma(1+2j)}
$$
starting from $\hu_{-1,*}(z)\equiv 0$ holds for all $j\geq 0$.
Applying the Nagumo norms of indices $(2j,r)$ and properties {\it (iv)}
and {\it (v)} of Proposition  \ref{propNagumo} we get
\begin{equation*}
\begin{array}{ccccl}
\displaystyle \frac{\Vert \hu_{j,*}(z)\Vert_{2j,r}}{\Gamma(1+2j)} &\leq &
\displaystyle \frac{\Vert
  \hf_{j,*}(z)\Vert_{2j,r}}{\Gamma(1+2j)}&+&\displaystyle \Vert
a(z)\Vert_{0,r}\frac{\Vert \hu''_{j-1,*}(z)\Vert_{2j,r}}{\Gamma(1+2j)}\\
\noalign{\medskip}
&\leq & ''&+&\displaystyle \Vert a(z)\Vert_{0,r}\, e^2\,  \frac{\Vert
  \hu_{j-1,*}(z)\Vert_{2j-2,r}}{\Gamma\big(1+(2j-2)\big)}\\
\end{array}
\end{equation*}
Denote $\displaystyle g_j=\frac{\Vert
  \hf_{j,*}(z)\Vert_{2j,r}}{\Gamma(1+2j)}$ and $\displaystyle \alpha= \Vert
a(z)\Vert_{0,r}\, e^2$ and consider the numerical
sequence $$\displaystyle \left\lbrace \begin{array}{l}
v_{-1}=0\\
\noalign{\medskip}
v_j=g_j+\alpha \,v_{j-1} \text{ for all } j\geq 0.
\end{array}\right.$$
By construction, $\displaystyle  \frac{\Vert
  \hu_{j,*}(z)\Vert_{2j,r}}{\Gamma(1+2j)} \leq v_j$ for all $j\geq 0$.\\
Let us bound $v_j$ as follows. By assumption,
$\displaystyle 0\leq g_j\leq \frac{CK^j\Gamma(1+2j)}
{\Gamma(1+2j)}r^{2j}
=C(Kr^2)^j$ for all $j$ and the series $
g(X)=\sum_{j\geq 0} g_j X^j$ is convergent. Due to the
recurrence relation defining the $v_j$'s the series
$v(X)=\sum_{j\geq 0} v_jX^j$ satisfy $(1-\alpha
X)v(X)=g(X)$. It is then convergent and there exist
constants $C'>0, K'>0$ such that $v_j\leq C' \, {K'}^j$ for
all $j$. Hence,
$$
\Vert \hu_{j,*}(z)\Vert_{2j,r} \leq C'\, {K'}^j \Gamma(1+2j)\quad
\text{for all } j\geq 0.
$$
We deduce a similar estimate on the sup-norm by shrinking the domain
$D_r$. Indeed, let $0<r'<r$. For all $j\geq 0$ and $z\in D_{r'}$,
$$\begin{array}{ccl}
\vert \hu_{j,*}(z)\vert &=&\displaystyle \Big\vert \hu_{j,*}(z)d_r(z)^{2j}
\frac{1}{d_r(z)^{2j}}\Big\vert\\
\noalign{\medskip}
&\leq&\displaystyle \frac{1}{(r-r')^{2j}}\big\vert
\hu_{j,*}(z)d_r(z)^{2j}\big\vert\\
\end{array}
$$
Hence,
$$\begin{array}{ccl} \displaystyle
\sup_{z\in D_{r'}} \vert \hu_{j,*}(z)\vert &\leq&\displaystyle
 \frac{1}{(r-r')^{2j}}\, \Vert \hu_{j,*}\Vert_{2j,r}\\
\noalign{\medskip}
&\leq&\displaystyle C'\,\Big(\frac{K'}{(r-r')^2}\Big)^j \Gamma(1+2j)
\end{array}
$$
\qed
\section{1-summability}
Still considering $t$ as the variable and $z$ as
a parameter, one  extends the classical notions of
summability to families parameterized by $z$ in requiring
similar conditions, the estimates  being however
uniform with respect to the parameter $z$. For a general study of series with coefficients in a Banach space we refer to \cite{Bal00}. Among the
many equivalent definitions of 1-summability in a given
direction $\arg t=\theta$ at $t=0$ we choose here a generalization of Ramis
definition which states that a series $\hf$ is
1-summable in the direction $\theta$ if there exists a
holomorphic function $f$ which is 1-Gevrey asymptotic  to
$\hf$ on an open sector $\Sigma_{\theta,>\pi}$  bisected by $\theta$
with opening larger than
$\pi$ (\cf \cite{R80} D\'ef 3.1).
 There are various equivalent ways of expressing the
1-Gevrey asymptoticity. We choose to extend the one which
sets conditions on the successive derivatives of $f$  (see
\cite{Mal95} p. 171 or \cite{R80} Thm 2.4, for instance).

\begin{defi}[1-summability]\label{1sum}
A series $\hu(t,z)\in\Ot$ is {\rm 1-summable in the direction
$\arg t=\theta$} if there exist a sector $\Sigma_{\theta,>\pi}$, a
radius $0<r\leq\rho$ and a function $u(t,z)$ called {\rm 1-sum of $\hu(t,z)$ in the direction $\theta$ }such that
\begin{enumerate}
\item $u$ is defined and holomorphic on
$\Sigma_{\theta,>\pi}\times D_r$;
\item For any $z\in D_r$ the map
$t\mapsto u(t,z)$ has
$\displaystyle
\hu(t,z)=\sum_{j\geq 0}\frac{t^j}{j!} \, \hu_{j,*}(z)$ as
Taylor series at 0 on $\Sigma_{\theta,>\pi}$;
\item For any proper\footnote{In this context a subsector
$\Sigma$ of a sector
$\Sigma' $ is said a proper subsector and one denotes
$\Sigma\subset\subset\Sigma'$ if its closure in
$\C$ is contained in $\Sigma'\cup \{0\}$.} subsector
$\Sigma\subset\subset \Sigma_{\theta,>\pi}$ there exist constants
$C>0, K>0$  such that for all $\ell\geq 0$, all $t\in \Sigma
$ and $z\in D_r$
$$
\begin{tabular}
{| p{5.8cm}  |}
\hline  $\displaystyle
\quad \big\vert \partial_t^\ell\, u(t,z)\big\vert \leq CK^\ell \Gamma(1+2\ell)
$
\tvi{20}{17}\\ \hline
\end{tabular} \ .
$$
\end{enumerate}
\end{defi}
We denote $\Ots_{1,\theta}$ the subset of $\Ot$ made of
all 1-summable series in the direction $\arg t=\theta$.
Actually, $\Ots_{1,\theta}$ is included in $\Ot_1$.
\medskip

For any fixed $z\in D_r$, the 1-summabilty of the series $\hu(t,z)$ is  the classical 1-summability and Watson Lemma implies the unicity of its 1-sum, if any.

\begin{pro}\label{diffalsum}
$\big( \Ots_{1,\theta}, \partial_t,\partial_z\big)$ is a differential
$\C$-algebra stable under $\partial_t^{-1}$ and $\partial_z^{-1}$.
\end{pro}
\proof Let $\hu(t,z)$ and $\hv(t,z)$ be two 1-summable series in
direction $\theta$. In Definition \ref{1sum} we can choose the same
constants $r, C, K$ both for $\hu$ and $\hv$. The product
$w(t,z)=u(t,z)v(t,z)$ satisfies conditions {\it 1} and {\it 2} of
Definition \ref{1sum}.
Moreover,
\begin{equation*}
\begin{array}{ccl}
\big\vert \partial_t^\ell w(t,z)\big\vert &=&\displaystyle
\Big\vert \sum_{p=0}^\ell \begin{pmatrix} \ell\\ p\\ \end{pmatrix}
\partial_t^p u(t,z)\partial_t^{\ell-p} v(t,z)\Big\vert\\
&\leq&\displaystyle  C^2\, K^\ell\, \Gamma(1+2\ell)\,\Bigg\vert
 \sum_{p=0}^\ell \frac{\Gamma(1+\ell)}{\Gamma(1+2\ell)}
 \frac{\Gamma(1+2p)}{\Gamma(1+p)}
 \frac{\Gamma\big(1+2(\ell-p)\big)}
{\Gamma\big(1+(\ell-p)\big)}\Bigg\vert\\
\noalign{\bigskip}
 &\leq&C^2 \,K^\ell \,(\ell+1)\,\Gamma(1+2\ell)\\
\noalign{\bigskip}
 &\leq & C'\,{K'}^\ell\, \Gamma(1+2\ell)\quad \text{for adequate } C', K'>0.\\
\end{array}
\end{equation*}
This proves condition {\it 3} of Definition \ref{1sum} for $w(t,z)$,
that is, stability of $\Ots_{1,\theta}$ under multiplication.\\
Stability under $\partial_t$, $\partial_t^{-1}$ or $\partial_z^{-1}$
is straightforward. Stability under $\partial_z$ is obtained using the
Integral Cauchy Formula on a disc $D_{r'}$ with $r'<r$.\qed
\bigskip

We may notice that the 1-sum $u(t,z)$ of a 1-summable series
$\hu(t,z)\in \Ots_{1,\theta}$ may be analytic with respect to $z$ on a
disc  $D_r$ smaller than the common disc $D_\rho$ of analyticity of
the coefficients $\hu_{j,*}(z)$ of $\displaystyle \hu(t,z)=\sum_{j\geq
  0}\frac{t^j}{j!}\hu_{j,*}(z)$. With respect to $t$, the 1-sum $u(t,z)$
is analytic on a sector supposedly open and containing a closed sector
$\overline{\Sigma}_{\theta,\pi}$ bisected by $\theta$ with opening
$\pi$; there is no control on the angular opening except that it must
be larger than $\pi$ and no control on the radius of this sector
except that it must be positive.
Thus, the 1-sum $u(t,z)$ is well defined as a section of the sheaf of
analytic functions in $(t,z)$ on a germ of closed sector of opening
$\pi$ (\ie a closed interval $\overline{I}_{\theta,\pi}$ of length
$\pi$ on the circle $S^1$ of directions issuing from 0, \cf \cite{MalR92} 1.1 or
\cite{L-R94} I.2) times $\{0\}\subset
\C_z$. We denote $\sO_{\overline{I}_{\theta,\pi}\times \{ 0\}}$ the
space of such sections.

\begin{cor}\label{cordiffalg}
The operator of 1-summation
$$
\sS :\ \left\lbrace\begin{array}{ccl} \Ots_{1,\theta}
    &\longrightarrow& \sO_{\overline{I}_{\theta,\pi}\times \{ 0\}}\\
\noalign{\medskip}
\hu(t,z) &\mapsto & u(t,z)\\
\end{array}\right.
$$
is a homomorphism of differential $\C$-algebras for the derivations
$\partial_t$ and $\partial_z$ and it commutes with $\partial_t^{-1}$
and $\partial_z^{-1}$.
\end{cor}

%%%%%%%%%%%%%%%%%%%%%%%%%%%%%%%%%%%%%%%
\begin{thm}\quad\\
Let a direction $\arg t=\theta$ issuing from 0 and a series
$\hf(t,z)\in\Ot$ be given.\\
Recall $D=1-a(z)\partial_t^{-1}\partial_z^2$ and assume that either
$a(0)\neq 0$ or $a(0)=0$ and $a'(0)\neq 0$.

Then, the unique solution $\hu(t,z)$ of $D\hu=\hf$ in $\Ot$ is
1-summable in the direction $\theta$ if and only if $\hu_{*,0}(t),
\hu_{*,1}(t)$ and $\hf(t,z)$ are 1-summable in the direction $\theta$.

Moreover, the 1-sum $u(t,z)$, if any, satisfies equation (\ref{Dheat}) in which $\hf(t,z)$ is replaced by the 1-sum $f(t,z)$ of $\hf(t,z)$ in direction $\theta$.
\end{thm}
%%%%%%%%%%%%%%%%%%%%%%%%%%%%%%%%%%%%%%
\proof  We first place ourselves in the case $a(0)\neq 0$.

Denote $\displaystyle a(z)=\sum_{n\geq 0}a_nz^n$.\\ As a preliminary
remark we notice that, by identification
of equal powers of $z$ in Equation
$$
\big( 1-a(z)\,\partial_t^{-1}\partial_z^2 \big)
\sum_{n\geq 0}\hu_{*,n}(t)\,\frac{z^n}{n!}=\sum_{n\geq 0}
\hf_{*,n}(t) \,\frac{z^n}{n!},\leqno{(\ref{Dheat})}
$$
we get
$$
\left\lbrace\begin{array}{l}
\hu_{*,0}(t)-a_0\, \partial_t^{-1}\hu_{*,2}(t) =\hf_{*,0}(t)\\
\noalign{\smallskip}
\hu_{*,1}(t)-a_1\, \partial_t^{-1}\hu_{*,2}(t)-a_0\,\partial_t^{-1}
\hu_{*,3}(t) =\hf_{*,1}(t)\\
\noalign{\smallskip}
\text{and so on}\dots
\end{array}\right.
$$
so that each $\hu_{*,n}(t)$ is uniquely and linearly
determined from
$\displaystyle \hu_{*,0}(t),\ \hu_{*,1}(t)$ and $\hf(t,z)$.
\begin{itemize}
\item The condition is necessary by Proposition
\ref{diffalsum}. Indeed, if $\hu$ is 1-summable then so
are $\displaystyle \hu_{*,0}(t)=\hu(t,0),\ \hu_{*,1}(t)=
\frac{1}{z}\big(\hu(t,z)-\hu_{*,0}(t)\big)\Big\vert_{z=0}$ and
$\hf=Du$.
\item Prove that the condition is sufficient. Assume that
$\displaystyle \hu_{*,0}(t),\ \hu_{*,1}(t)$ and $\hf(t,z)$ are
1-summable in direction $\theta$.

Set $\hu(t,z)=\hu_{*,0}(t)+z\,\hu_{*,1}(t)+
\partial_z^{-2}\hv(t,z)$ and
$\hw=\partial_t^{-1}\,\hv$.\\
With these notations Equation (\ref{Dheat}) becomes
\begin{equation}\label{fpm}
\Big(1-\frac{1}{a(z)}\partial_t\partial_z^{-2}
\Big)\,\hw(t,z)=\hg(t,z) \quad \text{where }\
\hg=\frac{1}{a(z)}(\hu_{*,0}+z\hu_{*,1}-\hf)
\end{equation}
and it suffices to prove that $\hw$ is 1-summable in direction
$\theta$ when $\hg$ is.
To this end, we proceed through a fixed point method as
follows.

Setting $\displaystyle \hw(t,z)=\sum_{p\geq 0} \hw_p(t,z)$
Equation  (\ref{fpm}) reads
$$
\begin{array}{l}\displaystyle
\ \ \hw_0 -\frac{1}{a(z)} \partial_t \partial_z^{-2}\, \hw_0 \quad = \hg\\
\noalign{\medskip}
\displaystyle +\hw_1 -\frac{1}{a(z)} \partial_t \partial_z^{-2}\, \hw_1\\
\noalign{\smallskip}
+\cdots\\
\displaystyle  +\hw_p-\frac{1}{a(z)}\partial_t \partial_z^{-2}\, \hw_p\\
\noalign{\smallskip}
+\cdots\\
\end{array}
$$
and we choose the solution given by the system
\begin{equation}\label{wp}
\left\lbrace\begin{array}{ccl}
\hw_0&=&\hg\\
\noalign{\smallskip}
\displaystyle \hw_1&=&\displaystyle
\frac{1}{a(z)}\partial_t \partial_z^{-2}\, \hw_0\\
\dots&&\\
\displaystyle \hw_p&=&\displaystyle
\frac{1}{a(z)}\partial_t \partial_z^{-2}\, \hw_{p-1}\\
\dots&&\\
\end{array}\right.
\end{equation}
We can check that, for all $p\geq 0$, the formal series
$\hw_p(t,z)$ are of order  $O(z^{2p})$ in $z$ and consequently,
 the series $\hw(t,z)=\sum_{p\geq 0}\hw_p(t,z)$ itself
makes sense as a formal series in $t$ and $z$.

Let $w_0(t,z)$ denote the  1-sum of $\hw_0=\hg$ in direction $\theta$ and
for all $p>0$, let $w_p(t,z)$ be determined as the solution of System
(\ref{wp}) in which all $\hw_p$ are replaced by $w_p$. All  $w_p$
are defined on a common domain $\Sigma_{\theta,>\pi}\times
D_{\rho'}$.\\
 We are willing to prove that the series $\displaystyle \sum_{p\geq 0}
w_p(t,z)$ is convergent with sum $w(t,z)$, the 1-sum of $\hw(t,z)$ in direction $\theta$.
 \\

The 1-summability of
$\hw_0$ implies that there exists $0<r'<\rho'$ and, for any subsector
$\Sigma\subset\subset
\Sigma_{\theta,>\pi}$, there exist constants $C'>0$, \, $K'>0$
such that for all $\ell\geq 0$ and $(t,z)\in
\Sigma\times D_{r'}$,
$$
\big \vert \partial_t^\ell w_0(t,z)\big\vert \leq C'{K'}^\ell\;
\Gamma(1+2\ell).
$$
Denote $\displaystyle B=\max_{z\in D_r}\Big\vert
\frac{1}{a(z)}\Big\vert$\\
From $\displaystyle w_1=\frac{1}{a(z)} \,
\partial_t\partial_z^{-2} w_0$ we deduce that
$$
\begin{array}{ccl}
\big\vert \partial_t^\ell
w_1\big\vert&=&\displaystyle\Big\vert\frac{1}{a(z)}\,
\partial_t^{\ell+1}\partial_z^{-2}\, w_0\Big\vert\\
\noalign{\medskip}
&\leq&\displaystyle B\max_{z\in
D_r}\big\vert\partial_t^{\ell+1}w_0\big\vert \frac{\vert
z\vert2}{2!}\\
\noalign{\medskip}
&\leq&\displaystyle
C'\,{K'}^{\ell+1}\,\Gamma\big(1+2(\ell+1)\big)\frac{B\vert
z\vert^2}{2!}\\
\end{array}
$$
and, by recursion, that
\begin{equation}\label{dlp}
\big\vert\partial_t^\ell w_p(t,z)\big\vert\leq C'{K'}^{\ell+p}\,
\Gamma\big(1+2(\ell+p)\big)\frac{(B\vert z\vert^2)^p}{(2p)!}
\quad \text{for all } p\geq 0.
\end{equation}
This implies
$$
\begin{array}{ccl}\displaystyle
\sum_{p\geq 0}\big\vert \partial_t^\ell w_p(t,z)\big\vert&\leq&\displaystyle
 C'\,{K'}^\ell\,
\Gamma(1+2\ell)\sum_{p\geq 0}
\left(\begin{matrix}2\ell+2p\\2p\\ \end{matrix}\right)
\big(K'\,B\,\vert z\vert^2\big)^p\\
&\leq&\displaystyle C'\,(4K')^\ell\,
\Gamma(1+2\ell)\sum_{p\geq 0} \big(4K'B\vert
z\vert^2\big)^p\\
&& \displaystyle \text{since }
\left(\begin{matrix}2\ell+2p\\2p\\ \end{matrix}\right)
\leq \sum_{k=0}^{2\ell+2p}\left(\begin{matrix}2\ell+2p\\k\\
\end{matrix}\right)=2^{2\ell+2p}. \\
\end{array}
$$
Denote $L=4K'Br^2$ and choose $r$ so small that $L<1$.\\
Denote $ C=C' \sum_{p\geq
0}L^p<\infty$ and $K=4K'$.\\
Then,
\begin{equation}\label{majorsigma}
\sum_{p\geq 0}\big\vert \partial_t^\ell w_p(t,z)\big\vert\leq
CK^\ell \Gamma(1+2\ell) \quad \text{ on } \Sigma\times D_r.
\end{equation}
In particular, for $\ell =0$, the series $\sum w_p(t,z)$ is normally convergent on $\Sigma\times D_r$. Consequently, its sum $w(t,z)$ exists and is analytic on $\Sigma\times D_r$. This proves point {\it 1} of Definition \ref{1sum} if we choose as sector $\Sigma \subset \Sigma_{\theta,>\pi}$ a sector bisected by $\theta$ with opening larger than $\pi$ .\\
For all $\ell\geq 1$, the series $\sum\partial_t^\ell w_p(t,z)$ is also normally convergent on $\Sigma\times D_r$ so that the series $\sum w_p(t,z)$ can be derivated termwise infinitely many times with respect to $t$ and the estimates (\ref{majorsigma}) imply
\begin{equation}\label{majorsum}
\big\vert \partial_t^\ell w(t,z)\big\vert\leq
CK^\ell \Gamma(1+2\ell) \quad \text{ on } \Sigma\times D_r
\end{equation}
which proves the condition {\it 3} of Definition \ref{1sum}.

Moreover, summing the Equations (\ref{wp}) for $w_p$ and the 1-sum $g(t,z)$ instead of $\hw_p$ and $\hg(t,z)$ we get $\displaystyle w(t,z)=g(t,z)+\frac{1}{a(z)}\sum_{p\geq 0}\partial_t\partial_z^{-2} w_p(t,z)=
g(t,z)+\frac{1}{a(z)}\partial_t\partial_z^{-2} w(t,z)$. Hence, $w(t,z)$ satisfies  Equation (\ref{fpm}) with right hand side $g(t,z)$ in place of $\hg(t,z)$.

Finally, the fact that all derivatives of $w(t,z)$ with respect to $t$ are
bounded on $\Sigma$ implies the existence of  $\displaystyle
\lim_{t\rightarrow 0\atop t\in\Sigma}\partial_t^\ell
w(t,z)$  for all $z\in D_{r}$ and hence the
existence of the Taylor series of $w$ at 0 on
$\Sigma$ for all $z\in D_{r}$.
Since $w(t,z)$ satisfies Equation (\ref{fpm}), so does its Taylor series.
Since Equation (\ref{fpm}) has a unique formal  solution $\hw(t,z)$, we can conclude that the
Taylor expansion of $w(t,z)$  is $\hw(t,z)$, which proves
part {\it 2} of Definition \ref{1sum}.

This achieves the proof of the 1-summability of $\hu(t,z)$ in direction $\theta$ in the case when $a(0)\neq 0$.
\item The fact that the 1-sum $u(t,z)$ of $\hu(t,z)$ in direction $\theta$ satisfies Equation
  (\ref{Dheat}) with right hand side the 1-sum $f(t,z)$
 of $\hf(t,z)$ instead of $\hf(t,z)$ is equivalent to the fact that $w(t,z)$ satisfies
 Equation (\ref{fpm}) with right hand side $g(t,z)$
instead of $\hg(t,z)$, which we proved above. It is also a consequence of Corollary \ref{cordiffalg}.
\end{itemize}

\vspace{.5cm}

In the case when
$a(0)=0$ and $a'(0)\neq 0$ the necessary condition again results from
Proposition \ref{diffalsum}. The fact that $u(t,z)$ satisfies Equation
(\ref{Dheat}) results from Corollary \ref{cordiffalg}. We sketch the
proof of the sufficient condition.

 Denote $a(z)=zA(z)$ with $A(0)\neq 0$.
\\
 In this
case, identification of equal powers of
$z$ shows that $\hu_{*,0}=\hf_{*,0}$ and that all
$\hu_{*,n}$ for
$n\geq 1$ are uniquely determined by $\hu_{*,1}$ and $\hf$.\\
We set again
$\hu(t,z)=\hu_{*,0}+z\hu_{*,1}+\partial_t\partial_z^{-2}\hw$ so
that $\hw$ satisfies the equation
\begin{equation}
\Big(1-\frac{1}{zA(z)}\partial_t\partial_z^{-2}\Big)\,
\hw(t,z)=\hg(t,z) \quad \text{where }
\hg=\frac{1}{A(z)}\Big(\hu_{*,1}+\frac{\hu_{*,0}-\hf}{z}\Big).
\end{equation}
Still, $\hg$ is a formal series, assumed to be 1-summable in
direction $\theta$ and we look for $\hw$ in the form
$\displaystyle \hw=\sum_{p\geq 0}\hw_p$ as previously. The
operator $\displaystyle \frac{1}{z}\partial_z^{-2}$ implies
that $\hw^p=O(z^p)$ instead of $O(z^{2p})$. If we denote
$\displaystyle B=\max_{z\in D_r}\frac{1}{\vert A(z)\vert}$,
then, for all $p$ and $\ell$,
$$
\big\vert \partial_t^\ell w_p\big\vert\leq
C'K'^{\ell+p}\Gamma\big(1+2(\ell+p)\big)\frac{(B\vert
z\vert)^p}{p!}
$$
and it follows that, for a
convenient choice of $r>0$,
$$
\big\vert \partial_t^\ell w(t,z)\big\vert\leq C\,
{K}^\ell\, \Gamma(1+2\ell)
$$
with $\displaystyle C=C' \sum_{p\geq
0}(4KBr)^p<\infty$  and
$K=4K'$.\qed

\vspace{.5cm}

 The case of a thermal diffusivity $a(z)=O(z^2)$ gives rise to the
conditions $\hu_{*,0}(t)=\hf_{*,0}(t)$ and $\hu_{*,1}(t)=\hf_{*,1}(t)$
and we could hope of similar necessary and sufficient conditions
which apply to the inhomogenuity $\hf(t,z)$ only. This is not the case
since the previous proof cannot be extended to that situation.
Indeed, the appearance of  $\displaystyle \frac{\partial_z^{-2}}{z^2}$
instead of $\displaystyle \partial_z^{-2}$ or $\displaystyle
\frac{\partial_z^{-2}}{z}$ implies that no power of
$z$ remains in the estimates (\ref{dlp}) and we cannot guaranty the
convergence of the estimate for $\partial_t^\ell w$.

 The counter-example below shows that even with $\hf(t,z)$ independent
 of $t$ and rational the 1-summability of $\hu(t,z)$ may fail.
\begin{cexa}\quad\end{cexa}

Consider the heat initial conditions problem (\ref{Dheat}) with   $\displaystyle \hf(t,z)=\sum_{n\geq 0} z^n=\frac{1}{1-z}$ and $a(z)\equiv 1$. The series $\hf(t,z)$  is independent of $t$ and is  convergent
 in $z$ near 0 with rational sum.
 The problem is equivalent to the heat initial conditions problem without internal heat generation
 \begin{equation}\left\lbrace\begin{array}{l}
 \partial_t\hu-z^2\,\partial_z^2 \hu=0\\
 \noalign{\medskip}
 \displaystyle \hu(0,z)=\sum_{n\geq 0} z^n\\
 \end{array}\right.
 \end{equation}
In this case, $\hu_{*,0}(t)=\hf_{*,0}(t)\equiv 1$, \
 $\hu_{*,1}(t)=\hf_{*,1}(t)\equiv 1$ and for all $n\geq 2$,
 $\hu_{*,n}(t)$ satisfies
 \begin{equation*}
\hu'_{*,n}(t)-n(n-1)\hu_{*,n}(t)=0\quad \text{and } \ \hu_{*,n}(0)=n!.
 \end{equation*}
Consequently, $\hu_{*,n}(t)=n!\, e^{n(n-1)t}$.
\bigskip

Suppose $\hu(t,z)$ is 1-summable in a direction $\theta$ with sum
$u(t,z)$. \\Then,
since $\displaystyle \hu_{*,n}(t)=\partial_z^n\hu(t,z)\Big\vert_{z=0}$ all
$\hu_{*,n}(t)$ are 1-summable in direction $\theta$ with sum
$u_{*,n}(t)=\partial_z^n u(t,z)\Big\vert_{z=0}$. The Integral Cauchy
Formula applied to $\partial_z^n u(t,z)$ at $z=0$ provides estimates
of the form
$$
\vert u_{*,n}(t)\vert =\Bigg\vert \frac{n!}{2\pi i} \int_{\vert
  \zeta\vert =R<r} \frac{u(t,\zeta)}{\zeta^{n+1}}d\zeta\Bigg\vert\leq
\frac{n!}{2\pi}\frac{C 2\pi R}{R^{n+1}}=C \, k^n\, n!
$$
on a sector bisected by $\theta$ with opening larger than $\pi$.
In our case, $\hu_{*,n}(t)=u_{*,n}(t)=n!\, e^{n(n-1)t}$. The functions
$e^{n(n-1)t}$ being  unbounded on any sector larger than a half plane
such estimates are impossible. Hence, $\hu(t,z)$ is 1-summable in no
direction.\qed

\section{Initial conditions}

We end this article with a discussion of how to apply the
above result and we  develop the cases when $a(z)=a\in\C^*$ or
$a(z)=bz, \ b\in\C^*$.
\bigskip

The formal series  $\hf(t,z)$ is a data of the  problem and
although its 1-summability may be not obvious we assume that
it is known. $\hf(t,z)$ is not itself the initial conditions but is
closely connected to (see Section 1).

The series $\hu_{*,0}(t)$ and $\hu_{*,1}(t)$ can, at least
theoretically, be computed in terms of $\hf(t,z)$ from the formula
$$\displaystyle \hu(t,z)=\sum_{k\geq 0}\big(
a\partial_t^{-1}\partial_z^2\big)^k \hf(t,z)$$ and an
explicit  computation can be achieved for simple $a(z)$
such as
$a(z)=a$ constant, $a(z)=bz$ ($b\in\C^*$) or $a(z)=a+bz$.
However, an explicit computation of
$\hu_{*,0}(t)$ and
$\hu_{*,1}(t)$ looks like hopeless  for a general $a(z)$.
 \protect\boldmath
\subsection{Case $a(z)=a\in\C^*$}\protect\unboldmath

When $a$ is a constant then the operators $a, \partial_t$ and
$\partial_z$ commute
and
$\big(a\partial_t^{-1}\partial_z^2\big)^k=a^k\partial_t^{-k}\partial_z^{2k}$.
 From the calculation of $ \hu(t,z)=\sum_{k\geq 0}\big(
a\partial_t^{-1}\partial_z^2\big)^k \hf(t,z)$ we obtain
\begin{equation}\label{u01a}
\left\lbrace\begin{array}{ccl}
\hu_{*,0}(t)&=&\displaystyle \sum_{k\geq
  0}\frac{t^k}{k!}\sum_{j+n=k}a^n \hf_{j,2n}\\
\noalign{\medskip}
\hu_{*,1}(t)&=&\displaystyle \sum_{k\geq
  0}\frac{t^k}{k!}\sum_{j+n=k}a^n \hf_{j,2n+1}\\
\end{array}\right.
\end{equation}
Our aim is to characterize the 1-summability of these two series as a
property of the inhomogenuity $\hf$.

\bigskip

\noindent $\bullet$ We start with the  case where $\displaystyle  \hf(t,z)=
\sum_{n\geq
  0}\hf_{0,n}\frac{z^n}{n!}$ is independent of $t$  which
corresponds to Problem (\ref{homheat}). For simplicity, we denote
$\hf(z)$.

The formul\ae\ (\ref{u01a}) become
\begin{equation}
\left\lbrace\begin{array}{ccl}
\hu_{*,0}(t)&=&\displaystyle \sum_{k\geq
  0}\frac{(at)^k}{k!} \hf_{0,2k}\\
\noalign{\medskip}
\hu_{*,1}(t)&=&\displaystyle \sum_{k\geq 0}\frac{(at)^k}{k!} \hf_{0,2k+1}\\
\end{array}\right.
\end{equation}

Define the 2-Laplace transform of $\hf(z)$ by
$\displaystyle \sL_z^{[2]}\hf(\gz)=\sum_{n\geq
  0}\hf_{0,n}\frac{\gz^n}{n!}\,\frac{n!}{[n/2]!}$ where $[n/2]$ stands
for the integer part of $n/2$. Then,
$$
\sL_z^{[2]}\hf\big((at)^{1/2}\big)=\hu_{*,0}(t)+(at)^{1/2} \hu_{*,1}(t).
$$
and we may state
\begin{pro}\label{ahomog}
Suppose $a(z)=a\in\C^*$ and $\hf(t,z)=\hf(z)$.\\
Then, the following three assertions are equivalent.
\begin{description}
\item{ (i) \ } \ $\hu_{*,0}(t)$ and $\hu_{*,1}(t)$ are 1-summable in
  direction $\theta$;
\item{(ii) } \  $\sL^{[2]}_z\hf(z)$ is 2-summable in the directions
  $\frac{1}{2}(\theta+\arg a)$ mod $\pi$;
\item{(iii)} \  $\hf(z)$ is analytic near 0 and it can be analytically
  continued to sectors neighbouring the directions
  $\frac{1}{2}(\theta+\arg a)$ mod $\pi$  with exponential growth of
  order 2 at infinity.
\end{description}
\end{pro}

Assertion {\it (iii)} with $a=1$ (hence $\arg a=0$) is how the conditions are
formulated in \cite{LMS99} and proved via direct Borel-Laplace
estimations.
Our method provides thus a new proof of this result.
\vspace{.5cm}

\bigskip

\noindent $\bullet$ Consider now the case of a general $\hf(t,z)$.\\
The interpretation of the 1-summability of $\hu_{*,0}(t)$ and
$\hu_{*,1}(t)$ becomes more involved and uses Borel and Laplace transforms
of $\hf(t,z)$ in both variables.
%Note first that the 2-Laplace transform
%$\sL^{[2]}_z\hf(z)$ above  is nothing else than a 1-Laplace followed
%by a 2-Borel transform with respect to $z$ (2-Borel {\it w.r.t.} $z$
%means 1-Borel {\it w.r.t.} $z^{1/2}$).
\\
We denote $\sL_z$ or $\sB_z$ and so on\dots the 1-Laplace or 1-Borel
transform  {\it w.r.t.} $z$ and so on\dots. These operators are
defined here by $\sL_zz^n=\zeta^n [n]!$ and $\sB_z=\sL_z^{-1}$ where $[n]$
denotes the integer part of $n$.\\
Consider $\displaystyle \tvi{20}{20}
\sL_t\sL_z\hf\big(\tau, (a\tau)^{1/2}\big)=\sum_{k\geq 0}
\tau^k\sum_{j+n=k}\hf_{j,2n} a^n + (a\tau)^{1/2}\sum_{k\geq 0}
\tau^k\sum_{j+n=k}\hf_{j,2n+1} a^n $ and\\ $\displaystyle
\sB_\tau\sL_t\sL_z\hf\big(\tau, (a\tau)^{1/2}\big)(t)=\sum_{k\geq 0}
\frac{t^k}{k!}\sum_{j+n=k}\hf_{j,2n} a^n + (at)^{1/2}\sum_{k\geq 0}
\frac{t^k}{k!}\sum_{j+n=k}\hf_{j,2n+1} a^n \tvi{0}{20} $ (the terms in
$\tau^k$ are divided by $k!$ and the terms in $\tau^{k+1/2}$ by $[k+1/2]!=k!$).\\
Denote $\hF(t)=\sB_\tau\sL_t\sL_z\hf\big(\tau,
(a\tau)^{1/2}\big)(t^2)$. Then,
$$
\hF(t^{1/2})=\hu_{*,0}(t)+(at)^{1/2} \hu_{*1}(t)
$$
and we may state:
\begin{pro}\label{ageneral}
Suppose $a(z)=a\in\C^*$ and $\hf(t,z)$ general.\\
Then, the series $\hu_{*,0}(t) $ and $\hu_{*,1}(t)$ are 1-summable in
direction $\theta$ if and only if the series $\hF$ associated with
$\hf$ as above is 2-summable in the directions $\theta/2$ mod $\pi$.
\end{pro}

The condition in Proposition \ref{ahomog} may be not easy to check but
seems reasonnable. In Proposition \ref{ageneral}, the link between
$\hf$ and $\hF$ is more complicated and the question remains of how to
check the 2-summability of $\hF$ in practice.

\protect\boldmath
\subsection{Case $a(z)=bz, \ b\in\C^*$}\protect\unboldmath

In this case,
$\big(a(z)\partial_t^{-1}\partial_z^2\big)^k=
b^k\partial_t^{-k}(z\partial_z^2)^k$ and
$$
(z\partial_z^2)^k\cdot\frac{z^n}{n!}=\left\lbrace\begin{array}{ll}
\displaystyle \frac{z^{n-k}}{(n-k)!}\frac{(n-1)!}{(n-k-1)!}& \text{ if } \ 0\leq
k<n\\
\noalign{\medskip}
0&  \text{ if } \ n\leq k.
\end{array}\right.
$$
 From the calculation of $ \hu(t,z)=\sum_{k\geq 0}\big(
bz\partial_t^{-1}\partial_z^2\big)^k \hf(t,z)$ we obtain
\begin{equation}\label{u01bz}
\left\lbrace\begin{array}{ccl}
\hu_{*,0}(t)&=&\displaystyle \sum_{j\geq
  0}\frac{t^j}{j!}\,\hf_{j,0}=\hf_{*,0}(t)\\
\noalign{\medskip}
\hu_{*,1}(t)&=&\displaystyle \sum_{j,k\geq
  0}\hf_{j,k+1}\,b^k\,\frac{t^{j+k}}{(j+k)!} k!\\
\end{array}\right.
\end{equation}
Since $\hu_{*,0}(t)=\hf_{*,0}(t)$ is 1-summable when so is $\hf(t,z)$,
our aim is now to characterize the 1-summability of the series
$\hu_{*,1}(t)$ as a
property of  $\hf$.
\bigskip

\noindent $\bullet$ Let us first again place  ourselves in the
situation of Problem (\ref{homheat})
where the inhomogenuity $ \hf(t,z)=
\sum_{n\geq 0}\hf_{0,n}\frac{z^n}{n!}$ is independent of $t$.

Formul\ae\ (\ref{u01bz}) become
\begin{equation}\label{u01abzhom}
\left\lbrace\begin{array}{ccl}
\hu_{*,0}(t)&=&=\hf_{0,0};\\
\noalign{\bigskip}
\hu_{*,1}(t)&=&\displaystyle \sum_{k\geq
  0}\hf_{0,k+1}\,b^k\, t^k.\\
\end{array}\right.
\end{equation}
Thus, $\sL_z \hf(bt)= \hf_{0,0}+ bt \hu_{*,1}(t)$ and we may state
\begin{pro}\label{bhomog}
Suppose $a(z)=bz, \ b\in\C^*$ and $\hf(t,z)=\hf(z)$.\\
Then, $\hu_{*,0}(t)$ is a constant and the following three assertions
are equivalent.
\begin{description}
\item{ (i) \ } \  $\hu_{*,1}(t)$ is 1-summable in
  direction $\theta$;
\item{(ii) } \  $\sL_z\hf(z)$ is 1-summable in the direction
  $\theta+\arg b$;
\item{(iii)} \  $\hf(z)$ is analytic near 0 and it can be analytically
  continued to a sector neighbouring the direction
  $\theta+\arg b$  with exponential growth of
  order 1 at infinity.
\end{description}
\end{pro}
\bigskip

\noindent $\bullet$  Consider the case of a general $\hf(t,z)$.\\
The Laplace transform of $\hf$ {\it w.r.t.} $z$ reads
$\ \sL_z\hf(t,z)=\hf_{*,0}(t)+z\sum_{j,n\geq 0}\frac{t^j}{j!}
\hf_{j,n+1} z^n$. Consider the series
$\hg(t,z)=\sL_t\sL_z\big\lbrack\frac{1}{z}\big(\sL_z\hf(t,z)-\hf_{*,0}(t)
\big) \big\rbrack $, We can check that  the Borel transform  of the
series $\hg(t,bt)$ is equal to $\hu_{*,1}(t)$ and we may state:
\begin{pro}\label{bzgeneral}
Suppose $a(z)=bz, \ b\in\C^*$ and $\hf(t,z)$ general.\\
Then, the series  $\hu_{*,1}(t)$ is 1-summable in
direction $\theta$ if and only if the Borel transform of $\hg(t,bt)$ is
1-summable in direction $\theta$.
\end{pro}
The comment following Propositions \ref{ahomog} and
\ref{ageneral} keeps valid.

%%%%%%%%%%%%%%%%%%%%%%%%%%%%

\end{document}